\documentclass[11pt]{article}%
   
\usepackage{amsmath,enumerate}
\usepackage{amsfonts}
\usepackage{amssymb}

\newtheorem{theorem}{Theorem}[section]

\newtheorem{corollary}[theorem]{Corollary}

\newtheorem{problem}[theorem]{Problem}
\newtheorem{proposition}[theorem]{Proposition}

\newtheorem{question}[theorem]{Question}

\newenvironment{proof}[1][Proof]{\noindent\textbf{#1.} }
{\hfill \ \rule{0.5em}{0.5em}}

\begin{document}
\title{Regular saturated graphs and sum-free sets}
\author{ Craig Timmons\thanks{Department of Mathematics and Statistics, California State University Sacramento, U.S.A. E-mail: \texttt{craig.timmons@csus.edu}. Research is supported in part by Simons Foundation Grant \#359419.}}

\maketitle

\begin{abstract}
In a recent paper, Gerbner, Patk\'{o}s, Tuza and 
Vizer studied regular $F$-saturated graphs.  
One of the essential questions is given $F$, for which 
$n$ does a regular $n$-vertex $F$-saturated graph exist.    
They proved that for all sufficiently large $n$, there 
is a regular $K_3$-saturated graph with $n$ vertices.   
We extend this result to both $K_4$ and $K_5$ and prove 
some partial results for larger complete graphs. 
Using a variation of sum-free sets from additive combinatorics, we 
prove that for all $k \geq 2$, there is a regular $C_{2k+1}$-saturated with 
$n$ vertices for infinitely many $n$.   
Studying the sum-free sets that give rise to $C_{2k+1}$-saturated graphs
is an interesting problem on its own and we state an open problem in this direction.    
\end{abstract}


\section{Introduction}

A graph $G$ is \emph{$F$-free} if 
$G$ has no subgraph isomorphic to $F$.  Understanding $F$-free graphs
is one of the central aims of extremal graph theory.
An interesting and important family of $F$-free graphs are
those that are maximal with respect to adding 
edges.  Such graphs are called \emph{$F$-saturated}.  
That is, a graph $G$ is 
$F$-saturated if $G$ is $F$-free and 
making any pair of nonadjacent vertices an edge
creates a copy of $F$.  
Let 
$\textup{sat}(n, F)$ denote the minimum number of edges in an 
$n$-vertex $F$-saturated graph.  This minimum is called 
the \emph{saturation number of $F$}.  
A general result of 
K\'{a}szonyi and Tuza \cite{kt} is that $\textup{sat}(n , F) = O_F(n)$ for
any graph $F$. Therefore, when $F$ has at least one edge $\textup{sat}(n,F)$ is 
linear in $n$.  The famous Erd\H{o}s-Hajnal-Moon Theorem \cite{ehm}
determines $\textup{sat}(n , K_s)$ exactly for all
$s \geq 3$ and characterizes the extremal graphs.   
This result is also attributed to Zykov \cite{zykov}. 
Determining $\textup{sat}(n , F)$ for other graphs 
$F$, such as cycles, complete bipartite graphs, and unions of cliques, is still an active area of research.  
The excellent surveys of Faudree, Faudree, and Schmitt \cite{sat survey} and Pikhurko \cite{pik} 
are standard references.  

In efforts to further understand saturation numbers and the structure of saturated graphs, 
there has been a substantial amount of research on 
adding degree conditions to the saturated graphs.  
Imposing a maximum degree or minimum degree condition 
has been investigated by 
 Alon, Erd\H{o}s, Holzman, Krivelevich \cite{aehk},
Day \cite{Day},
Duffus, Hanson \cite{DH}, 
F\"uredi, Seress \cite{fs}, Hanson and Seyffarth \cite{hs} to name a few.  
For example, the smallest maximum degree of a $K_s$-saturated 
graph with $n$ vertices is $\Theta (n^{1/2})$ (see \cite{aehk}), but 
an asymptotic formula is not known even in the case that
$s = 3$.  

Recently, 
Gerbner, Patk\'{o}s, Tuza, and Vizer
\cite{Gerbner-et-al} studied the existence of 
regular $F$-saturated graphs and their number of edges.  
Let 
\[
\textup{rsat}(n, F)
\]
be the minimum number of edges in an $n$-vertex regular 
$F$-saturated graph \textbf{when such a graph exists}.    
It was noted in \cite{Gerbner-et-al} that the 
principle concern is 
for which $n$ and $F$ does $\textup{rsat}(n,F)$ exist.
Gerbner et al. \cite{Gerbner-et-al} 
proved that $\textup{rsat}(n , K_3)$ exists for all $n \geq n_0$,
and gave some partial results for $\textup{rsat}(n, K_4)$.  
The aim here is to add to the number of graphs $F$ for which 
$\textup{rsat}(n,F)$ exists for infinitely many $n$.  

\begin{theorem}\label{K4 and K5 theorem}
For all $n > 59$, there is a regular $K_4$-saturated graph with $n$ vertices, 
and a regular $K_5$-saturated graph with $n$ vertices.    
\end{theorem} 

Our proof of Theorem \ref{K4 and K5 theorem} gives quadratic 
upper bounds on $\textup{rsat}(n , K_s)$~($s \in \{4,5 \}$) which we do not believe are optimal.  
Also, we did not attempt to minimize the lower bound on $n$.  It is included 
in the statement of Theorem \ref{K4 and K5 theorem} to show 
that we do not need $n$ to be large in our proof.     
For bigger cliques, we use the constructions from the proof of Theorem 
\ref{K4 and K5 theorem} and an observation of Gerbner et al.\ \cite{Gerbner-et-al} to prove 
a partial result.  

\begin{theorem}\label{all large values}
For any integer $\delta \geq 1$, there is a $k_{ \delta}$ such that the following 
holds for all $k \geq k_{ \delta}$.

\smallskip
\noindent
(i) For each $r \in \{0,1,2 , \dots , 2 \delta \}$, 
there is a $K_{2 \delta +1}$-saturated graph with $3 \delta k + r$ 
vertices that is $[ (3 \delta - 2) k + ( r - 1) ]$-regular.

\smallskip
\noindent
(ii) For each $r \in \{ 1,2, \dots , 2 \delta + 1 \}$, there 
is a $K_{2 \delta + 2 }$-saturated graph with 
$(3 \delta + 2)k +r$ vertices that is 
$[ 3 \delta k + (r - 1) ]$-regular.  
\end{theorem}

As $\delta$ gets large, Theorem \ref{all large values} 
shows that for roughly $\frac{2}{3}$'s of the possible residue classes of $n$ 
modulo $3 \delta$, there is a regular $n$-vertex $K_{2 \delta +1}$-saturated graph.
A similar statement can be made for regular $K_{2 \delta +2}$-saturated graphs.  

Next we look at odd cycles.
Let $\Gamma$ be an abelian group and $S \subseteq \Gamma$.  For an integer $k \geq 1$,  
the \emph{$k$-fold sumset of $S$} is 
\[
kS = \underbrace{S+S + \dots +S}_k = \{ s_1 + s_2 + \dots + s_k : s_i \in S \}.
\]
We say that $S$ is \emph{sum-free} if 
$2S \cap S = \emptyset$.  Sum-free sets 
have been studied extensively and we refer the reader to 
the survey of Tao and Vu \cite{tao-vu}.
A sum-free set $S \subseteq \Gamma$ is \emph{complete} if
\[
S + S = \Gamma \backslash S,
\]
and is \emph{symmetric} if $S = - S := \{ -s : s \in S \}$.
It is known (Section 1.2 of \cite{Haviv-Levy}) that if 
$S \subseteq \Gamma$ is a symmetric complete sum-free set, \textbf{s.c.sf}.\ for short, 
then the Cayley graph $\textup{Cay}( \Gamma , S)$ is 
a $K_3$-saturated graph that is $|S|$-regular and has $| \Gamma |$ vertices.
Theorem 1.5 of Haviv and Levy \cite{Haviv-Levy} and Proposition 
1.3 of \cite{Gerbner-et-al} implies that for large enough $n$, 
\[
\textup{rsast}(n , K_3)  = \Theta ( n^{3/2} ).
\]
Determining the correct constant is an open problem, but the point is that 
there is a strong connection between s.c.sf.\ sets and regular $K_3$-saturated graphs.  
Thus, it is reasonable to look for a variation of s.c.sf.\ sets that could be used to construct other
regular $F$-saturated graphs.  This led us to the next definition.  For 
a set $S \subseteq \Gamma$, let 
\[
\mathcal{R}_k (S) = \{ s_1 + s_2 + \dots + s_k : s_i \in S 
~\mbox{and}~
s_i + s_{i+1} + \dots + s_j \neq 0 
~\mbox{for all}~
1 \leq i < j \leq k
 \}.
\]
The set $\mathcal{R}_k (S)$ is all sums in $kS$ that can be ordered so that 
no consecutive sub-sum is 0. 
From the additive combinatorics perspective, $kS$ is a more natural object. 
For our application to graphs, we need 
$\mathcal{R}_k (S)$ because $k$-sums with no consecutive sub-sum being 0 allows
us to find paths of length $k$ from one vertex to another.  

\begin{proposition}\label{construction}
Let $k \geq 2$ be an even integer and let $\Gamma$ be an 
abelian group with $| \Gamma | = n$.  If there is a symmetric subset $S \subseteq \Gamma$
such that 
\[
\mathcal{R}_k ( S) = \Gamma \backslash  ( S \cup \{ 0 \} ) ~~ \mbox{and} ~~ 0 \notin (k+1)S,
\]
then the Cayley graph $\textup{Cay}( \Gamma , S)$ is a 
$C_{k+1}$-saturated $|S|$-regular graph with $n$ vertices.
\end{proposition}
   
The next theorem gives an infinite family for which we can 
apply Proposition \ref{construction}.  

\begin{theorem}\label{super sum set}
Let $\alpha , k  \geq 1$ be integers and $n = 2 \alpha ( \alpha + 4)k + 2 \alpha + 5$.  The subset
\[
S = \{ \pm ( 2 \alpha \ell + 1) : 0 \leq \ell \leq k \} \subset \mathbb{Z}_n
\]
has the property that $\mathcal{R}_{2 \alpha + 2 } ( S ) = \mathbb{Z}_n \backslash ( S \cup \{ 0 \} )$ and
$0 \notin (2 \alpha + 3)S$.  
\end{theorem}

\begin{corollary}\label{odd cycle corollary}
For all positive integers $\alpha$ and $k$,
there is a $C_{ 2 \alpha + 3}$-saturated $2 ( k + 1)$-regular $n$-vertex graph where
 $n = 2 \alpha ( \alpha + 4) k  + 2 \alpha + 5$.
\end{corollary}

Our proof of Theorem \ref{super sum set} shows that $(2 \alpha + 2 ) S = \mathbb{Z}_n \backslash S$
and we believe this raises an interesting problem which has been discussed in the additive combinatorics 
literature.  A subset $S \subseteq \Gamma$ of an abelian group $\Gamma$ is called a \emph{$(k , \ell )$-sum-free}
set if $(k S ) \cap ( \ell S) = \emptyset$.  For more on these types of sets and related problems, 
see \cite{Bajnok2, Bajnok, Bajnok-Matzke, HP}.  A $(k , \ell)$-sum-free set $S$ is 
\emph{complete} if $kS \cup \ell S $ is a partition of $\Gamma$.   The following 
is then a question about symmetric complete $(k,1)$-sum-free sets.  

\begin{problem}
For which $n$ and $k$ is there a symmetric set $S \subset \mathbb{Z}_n$ such that 
\[
kS = \mathbb{Z}_n \backslash S ~?
\]
\end{problem}

The upper bound on $\textup{rsat}( 2 \alpha ( \alpha + 4) + 2 \alpha + 5 , C_{ 2 \alpha + 3} )$
implied by Corollary \ref{odd cycle corollary} is quadratic in $n$ with the constant depending on $\alpha$.
Even in the case of $C_5$, it would be interesting to determine if there is a sparse version of $S$ that 
would give a subquadratic bound on $\textup{rsat}(n , C_5)$ for infinitely many $n$.  As mentioned 
above, the work of Haviv and Levy \cite{Haviv-Levy} achieves this for $K_3$ and gives an upper 
bound of the correct order of magnitude.  

The remainder of this paper is organized as follows.  In Section \ref{section k3} we prove 
Proposition \ref{construction} and Theorem \ref{super sum set}.  These two results prove Corollary \ref{odd cycle corollary}.
In Section \ref{ks section} we prove Theorems \ref{K4 and K5 theorem} and \ref{all large values}.  
Some concluding remarks are made in Section \ref{conclusion}.

   
\section{Proof of Proposition \ref{construction} and Theorem \ref{super sum set}}\label{section k3}


\noindent
\begin{proof}[Proof of Proposition \ref{construction}]
Suppose $\Gamma$ is an abelian group, written additively, and $S \subseteq \Gamma$ is a symmetric 
subset with $\mathcal{R}_k (S) = \Gamma \backslash ( S \cup \{ 0 \} )$
and  $0 \notin (k+1) S$.
First we show that the Cayley graph $G := \textup{Cay}( \Gamma , S)$ is 
$C_{k+1}$-free.   
For contradiction, suppose that $g_1 g_2 \dots g_{k+1} g_1$ 
is a $C_{k+1}$ in $G$ so
\begin{center}
$g_2 - g_1 = s_1$, ~~$g_3 - g_2 = s_2$, $\dots$, $g_{k+1} - g_k = s_k$, ~and ~
$g_1 - g_{k+1} = s_{k+1}$
\end{center}
for some $s_i \in S$.  We then have 
\begin{eqnarray*}
0 & = & (g_2 - g_1) + (g_3 - g_2) + \dots + ( g_{k+1} - g_k ) + ( g_1 - g_{k+1} ) \\
& = & s_1 + s_2 + \dots + s_k + s_{k+1}
\end{eqnarray*}
which contradicts the fact that $0 \notin (k+1) S$.  

The next step is to show that any two non-adjacent vertices are joined by a path 
of length $k$.
Suppose that $g_1$ and $g_{k+1}$ are non-adjacent vertices in $G$
so that $g_{k+1} - g_1 \in \Gamma \backslash ( S \cup \{0 \} )$.  
By hypothesis,  
there is a sum $s_1 + \dots + s_k \in \mathcal{R}_k (S)$ with $s_1 + \dots + s_k = g_{k+1} - g_1$ and 
ordered so that no consecutive sub-sum is 0.  
For $2 \leq i \leq k$, consider the vertices 
\[
g_i := g_1 + s_1 + s_2 + \dots + s_{i-1}.
\]
Note that $g_{k+1} = g_1 + s_1 + s_2 + \dots + s_k$ so 
the above equation holds for $i = k+1$ as well.  
We now claim that each of $g_1,g_2, \dots , g_k$, and $g_{k+1}$ are all distinct.
For contradiction, suppose that $g_j = g_i$ for some $1 \leq i < j \leq k+1$
where we can assume that $(i,j) \neq (1,k+1)$ since $g_1 \neq g_{k+1}$.
We then have 
\[
0 = g_j - g_i = s_{j-1} + s_{j-2} + \dots + s_{i-1}
\]
contradicting the fact that the sum $s_1 + \dots + s_k$ is in 
$\mathcal{R}_k (S)$.  We conclude that $g_1,g_2, \dots , g_{k+1}$ are 
all distinct.  By definition, $g_{i+1} - g_i = s_i$ for $1 \leq i \leq k$ 
which implies $g_1g_2 \dots g_k g_{k+1}$ is a path of length $k$ 
from $g_1$ to $g_{k+1}$.  
\end{proof}

\bigskip

\begin{proof}[Proof of Theorem \ref{super sum set}]
Let $ \alpha , k \geq 1$ be integers and $n = 2 \alpha ( \alpha + 4) k  + 2 \alpha + 5$.  
The cases $\alpha = 1$ and $\alpha \geq 2$ behave differently.  Let us assume that $\alpha \geq 2$ 
and then after completing the proof in this case, we will address the $\alpha =1$ case.   
At certain steps in the proof we will view the elements of $\mathbb{Z}_n$ as 
integesr.  When we do so, we will always use  
the least residues $\{0,1,2, \dots ,n  - 1 \}$ and the only time we do the arithmetic 
in $\mathbb{Z}$ is when we have been careful enough to ensure that there is no 
reduction modulo $n$, i.e., ``wrap around."

Let $S^+ =  \{  2 \alpha \ell + 1 : 0 \leq \ell \leq k \}$ 
and $S^- = \{ n - ( 2 \alpha \ell + 1) : 0 \leq \ell \leq k \}$ so that $S = S^+ \cup S^-$.  

Define $T_1 = \{ s_1 + \dots + s_{ 2 \alpha + 2} : s_i \in S^+ \}$, and for $2 \leq i \leq \alpha$, 
\begin{eqnarray*}
T_i & : = & \{  \underbrace{-1 + (-1) + \dots + (-1)}_{i-1} + s_i + s_{i+1} + \dots + s_{2 \alpha + 2} : 
s_i, s_{i+1}, \dots , s_{2 \alpha + 2}  \in S^+ , s_i \neq 1 \} \\
& = & 
\{  2  \alpha + 4 - 2 i + 2 \alpha ( \ell_i + \ell_{i+1} + \dots + \ell_{2 \alpha + 2} ) :
0 \leq \ell_i , \ell_{i+1}  , \dots , \ell_{ 2 \alpha + 2} \leq k , \ell_i \neq 1 \} \\
& = & 
\{ 2  \alpha + 4 -2i + 2 \alpha m : 1 \leq m \leq (2 \alpha + 3 - i )k \}.
\end{eqnarray*}
The set $T_i$ is a subset of $(2 \alpha + 2) S$, but we need sums in 
$\mathcal{R}_{ 2 \alpha + 2 } ( S)$.
We will define a set $U_i$ that is a truncated version of $T_i$ where we only consider the sums in $T_i$ that 
do not ``wrap around".  We need a bit more when $i > 1$ because in this case, 
we do not want $-1 + ( -1) + \dots + (-1)$~($i-1$ terms) to cancel with 
$s_i + s_{i+1} +\dots + s_{2 \alpha + 2}$.  This is where the condition $s_i \neq 1$ comes in.  
The requirement $s_i \neq 1$ implies that $s_i \in \{ 2 \alpha + 1 , 4 \alpha + 1 , \dots , 2 \alpha k + 1 \}$ so
\[
\underbrace{ -1 + (-1) + \dots + (-1) }_{i-1} + s_i \geq 2 \alpha + 2  - i \geq \alpha + 2 > 0.
\]  
This means that as long as we keep the sum $s_i + s_{i+1} + \dots + s_{2 \alpha + 2}$ small enough so that 
$- (i-1)  + s_i + s_{i+1} + \dots + s_{2 \alpha + 2} < n$, then we will have a sum in $\mathcal{R}_{2 \alpha + 2}(S)$
because 
no consecutive sub-sum of $- (i-1)  + s_i + s_{i+1} + \dots + s_{2 \alpha + 2} $ will be 0 modulo $n$.  
Before defining the sets $U_i$, we state a simple inequality:  
\begin{equation}\label{simple inequality}
( \alpha + 4) k \leq ( 2 \alpha + 3 - i ) k  ~~~ \Longleftrightarrow ~~~ i \leq \alpha - 1.
\end{equation}

Let 
\[
U_1 = \{ 2 \alpha +2   + 2 \alpha m : 1 \leq m \leq ( \alpha + 4)k \}
\]
and note that by (\ref{simple inequality}), $U_1$ is a subset of $T_1$.  The upper bound on $m$ implies that the largest element of 
$U_1$ is 
\[
2 \alpha + 2 +2 \alpha ( \alpha + 4) = n - 3.
\]
Therefore, the set $U_1$ is all $x \in \mathbb{Z}_n$
for which $x \equiv 4 - 2i ( \textup{mod}~2 \alpha )$ 
with the exception of 2.  
Next, let $V_1 = n - U_1$ which is all $x \in \mathbb{Z}_n$ for which 
\[
x \equiv n - ( 2 \alpha + 2 + 2 \alpha m ) \equiv 5 - 2 \equiv 3 ( \textup{mod}~2 \alpha )
\]
with the exception of $n - 2$.  In other words, 
\[
U_1 = \{ 2 \alpha + 2 , 4 \alpha + 2 , 6 \alpha + 2 , \dots , n - 3 \}
=
\{ x \in \{0,1, \dots , n - 1 \} : x \equiv 2 ( \textup{mod}~2 \alpha ) \} \backslash \{ 2 \}
\]
and
\[
V_1 = \{ 3 , 2 \alpha + 3 , 4 \alpha + 3 , \dots , n - (2 \alpha + 2) \}
=
\{ x \in \{0,1, \dots , n - 1 \} : x \equiv 3 ( \textup{mod}~2 \alpha ) \} \backslash \{ n - 2 \}.  
\]

For $2 \leq i \leq \alpha - 1$, let 
\[
U_i = \{ 2 \alpha + 4 - 2i + 2 \alpha m : 1 \leq m \leq ( \alpha + 4 )k \}.
\]
By (\ref{simple inequality}), $U_i$ is a subset of $T_i$.  The elements of $U_i$ are all
$x \in \mathbb{Z}_n$ for which $x \equiv 4 - 2 i ( \textup{mod}~ 2 \alpha )$ with the 
exception of $2 \alpha + 4 - 2i$.  
Let $V_i =  n- U_i$.  Then $V_i$ is all $x \in \mathbb{Z}_n$ for which 
\[
x \equiv n - ( 4 - 2i ) \equiv 5 - 4 + 2i \equiv 1 + 2i ( \textup{mod}~2 \alpha )
\]
with the exception of $n - (2 \alpha + 4 - 2i )$.  Equivalently, 
\[
U_i  = \{ x \in \{0,1, \dots , n - 1 \} : x \equiv 4 - 2i ( \textup{mod}~2 \alpha ) \} \backslash \{ 2 \alpha + 4 - 2i \}
\]
and
\[
V_i = \{ x \in \{ 0 , 1, \dots , n - 1 \} : x \equiv 1 + 2i ( \textup{mod}~2 \alpha ) \} \backslash \{ n  - ( 2 \alpha + 4 - 2i ) \}.
\]

Let $U_{ \alpha } = \{ 4 + 2 \alpha m : 1 \leq m ( \alpha + 3) k \}$.  This 
is a subset of $T_{ \alpha}$ since $( \alpha + 3)k \leq (2 \alpha + 3 - \alpha ) k$ and in fact, $U_{ \alpha} = T_{ \alpha}$ 
in this case.  The elements of $U_{ \alpha}$ are all $x \in \mathbb{Z}_n$ for which 
$x \equiv 4 ( \textup{mod}~2 \alpha)$ with the exeption of 4 and all elements of 
\[
S^-= \{ n - ( 2 \alpha k + 1) ,  n  - ( 2 \alpha ( k - 1) + 1) , \dots , n - (2 \alpha + 1) \}.
\]
Indeed, the largest element of $U_{ \alpha }$ is 
\[
4 + 2 \alpha ( \alpha + 3) k = n - ( 2 \alpha k + 1) - 2 \alpha.
\]
As before, let $V_{ \alpha } = n - U_{ \alpha}$ which is all $x \in \mathbb{Z}_n$ such that 
\[
x \equiv n- ( 4 - 2 \alpha ) \equiv 5 - 4 \equiv 1 ( \textup{mod}~2 \alpha)
\]
which the exception of $n - 4$ and all elements in
\[
S^+ = \{ 1 , 2 \alpha + 1 , \dots , 2 \alpha ( k - 1) + 1 , 2 \alpha k + 1 \}.
\]
By construction, all elements in each $U_i$ and $V_i$ are sums in $\mathcal{R}_{2 \alpha + 2 } (S)$. 
An elementary number theory argument shows that 
as $i$ ranges over $\{2, \dots , \alpha \}$, the values $4 - 2i ( \textup{mod}~2 \alpha )$ 
and $1 + 2i ( \textup{mod}~2 \alpha )$ cover all of the residues modulo $2 \alpha$ except for 2 and 3.  These
are covered by $U_1$ and $V_1$.  Thus,
\[
\bigcup_{i = 1}^{ \alpha } ( U_i \cup V_i)
= 
\mathbb{Z}_n \backslash (  S^+ \cup S^- \cup \{ 0 \} \cup \{2 , 4, \dots , 2 \alpha \} \cup \{ n - 2 , n -4 , \dots , n - 2 \alpha \}).
\]

Next we show that each element in $\{2,4, \dots , 2 \alpha \} \cup \{ n -2 , n -4 , \dots , n -2 \alpha \}$ is 
a sum in $\mathcal{R}_{ 2 \alpha + 2} (S)$.   
Since $k \geq 1$, the set $S$ contains $\{ \pm 1 , \pm ( 2 \alpha + 1 ) \}$.  Let 
$1 \leq j \leq \alpha$.  We can write 
\[
2j \equiv \underbrace{ - ( 2 \alpha + 1) - \dots - ( 2 \alpha + 1) }_{ \alpha + 1 -j } + 
\underbrace{1 + 1 + \dots + 1}_{2j } + \underbrace{ ( 2\alpha + 1) + \dots + ( 2 \alpha + 1) }_{ \alpha + 1 -j } ( \textup{mod}~n).
\]
Any consecutive sub-sum that is 0 must 
use all of the 1's since $2 \alpha + 1 > 2j$ for any $1 \leq j \leq \alpha$.  For contradiction, suppose that 
\[
- \theta_1 ( 2 \alpha + 1) +1 + \dots + 1 + \theta_2 ( 2 \alpha + 1) \equiv 0 ( \textup{mod}~n)
\]
where $0 \leq  \theta_1 , \theta_2 \leq \alpha -1 - j$.  If $\theta_1 = \theta_2$, then by rearranging the sum and cancelling 
we obtain $2j \equiv 0 ( \textup{mod}~n)$ which is a contradiction since $1 <  2j \leq 2 \alpha < n$.  
If $\theta_2 > \theta_1$, then we may cancel to get $0 \equiv 2j + ( \theta_2 - \theta_1) (2 \alpha + 1)$.
This is a contradiction since the right hand side is at most $2 + \alpha ( 2 \alpha  +1 ) < n$.
Finally, if $\theta_1 > \theta_2$ then we have $0 \equiv - ( \theta_1 - \theta_2)( 2 \alpha + 1 ) + 2j$.
This is again a contradiction since $0 > 2j - ( \theta_1 - \theta_2 ) (2 \alpha + 1) > -n$.  
The conclusion is that $2,4,6, \dots , 2 \alpha$ are all sums in $\mathcal{R}_{2 \alpha + 2}(S)$ and since 
$S$ is symmetric, we can say the same for $n-2 , n- 4, \dots , n- 2 \alpha$.
Therefore, $\mathbb{Z}_n \backslash ( S \cup \{ 0 \} ) \subseteq \mathcal{R}_{2 \alpha + 2} (S)$.  

We will complete the proof by showing that $0 \notin (2 \alpha + 3) S$.
This implies that $S \cap (2 \alpha + 2) S = \emptyset$ which in turn, implies 
$S \cap \mathcal{R}_{ 2 \alpha + 2} (S) = \emptyset$  since $\mathcal{R}_{ 2 \alpha + 2 } (S) \subseteq (2 \alpha + 2)S$.  
Suppose that  
\begin{equation}\label{congruence}
s_1 + \dots + s_{2 \alpha + 3} \equiv 0 ( \textup{mod}~n)
\end{equation}
for some $s_i \in S$.  
In $\mathbb{Z}$, 
\[
| s_1 + \dots + s_{ 2 \alpha + 3} | \leq ( 2 \alpha + 3) ( 2 \alpha k + 1) < 2n.
\]
This inequality with (\ref{congruence}) gives $s_1 + \dots + s_{2 \alpha + 3} \in \{ 0 , \pm n \}$. 
If $s_1 + \dots + s_{2 \alpha + 3} = 0$, then taking this equation modulo $2$ gives a contradiction since
each $s_i$ is odd.  Suppose now that $s_1 + \dots + s_{ 2 \alpha + 3} = n$.  Taking 
this equation modulo $2 \alpha$ leads to 
\[
2 \alpha + 3 \equiv 5 ( \textup{mod}~2 \alpha ) ~~~ \Longrightarrow ~~~ 3 \equiv 5 ( \textup{mod}~2 \alpha ).
\]
This last congruence implies $\alpha = 1$, but recall that we have assumed $\alpha \geq 2$.  This
is a contradiction.  The final possibility is $s_1 + \dots + s_{2 \alpha + 3} = -n $.
Multiplying this equation through by $-1$ and using the fact that $S$ is symmetric, we obtain 
$s_1 ' + \dots + s_{2 \alpha + 3}' = n$ where $s_i' \in S$.  As we have shown already, this 
is a contradiction.  

This completes the proof of the Theorem \ref{super sum set} when $\alpha \geq 2$. 
If $\alpha = 1$, then $n = 10k + 7$.  
We let $U_1$ be all sums of the form $s_1 + s_2 + s_3 + s_4$ 
with $s_i \in \{ 1 , 3, 5, \dots , 2k + 1 \}$ so $U_1 = \{ 4,6, 8 , \dots , 8k + 4 \}$.  Define 
 $V_1 = n - U_1$ so $V_1 = \{ 2k + 3 , 2k + 5 , \dots , 10k + 1 , 10k + 3 \}$. 
The elements of $U_1$ and $V_1$ are sums in $\mathcal{R}_4 (S)$.  
Since $k \geq 1$, we can write $2 = -3 + 1  + 1 + 3$ and $-2 = 3 - 1 - 1 -3$
which are also sums in $\mathcal{R}_4 (S)$.
Hence, $\mathbb{Z}_n \backslash ( S \cup \{ 0 \} ) \subseteq \mathcal{R}_4 ( S)$.  
A simple modification of the proof given above shows that $0 \notin 5S$.  
This completes the proof of the case when $\alpha = 1$. 
\end{proof}


\section{Constructing regular $K_s$-saturated graphs}\label{ks section}


In this section we prove Theorems \ref{K4 and K5 theorem} and \ref{all large values}.
The idea of the proof is to construct some families of regular $K_s$-saturated with $s \in \{3,4 \}$, and then use 
a Proposition of Gerbner et al.\ \cite{Gerbner-et-al}.  

\begin{proposition}[Gerbner, Patk\'{o}s, Tuza, Vizer]\label{useful tool}
 If $G$ is a $d_G$-regular $K_s$-saturated graph with $n_G$ vertices, and 
$H$ is a $d_H$-regular $K_t$-saturated graph with $n_H$ vertices, then the join $G+H$
is a $K_{s+t-1}$-saturated graph.  
If also $d_H + n_G = d_G + n_H $, then $G+H$ is $(d_H  + n_G)$-regular.  
\end{proposition}

Using a theorem of Haviv and Levy \cite{Haviv-Levy} on s.c.sf.\ sets,
we prove the an existence theorem on regular $K_3$-saturated graphs.
It is very important to mention that in \cite{Gerbner-et-al} a different construction 
of regular $K_3$-saturated graphs on $n \geq n_0$ vertices is given.  

\begin{theorem}\label{all K3}
Let $n$ be an integer with $n \geq 20$ or 
with $n \in \{5,8,10,11,13,14,16,17 \}$ and write $n = 3k + r$ where $r \in \{0,1,2 \}$.
For each such $n$, there is an $n$-vertex $d$-regular $K_3$-saturated graph
where $d = k + (r-1)$.
\end{theorem}
\begin{proof}
If $n = 3k+2$, then $S = \{ k+1,k+2, \dots , 2k +1 \}$
is a s.c.sf. set in $\mathbb{Z}_n$ (see \cite{Haviv-Levy}).  
For $k \geq 1$, the graph 
$\textup{Cay}( \mathbb{Z}_{3k+2} , S )$ will be $K_3$-saturated, have 
$3k+2$ vertices, and will be ($k+1$)-regular.  Similarly,
if $n = 3k+1$ with $k \geq 4$, then
\[
S = \{  k , 2k +1 \} \cup \{ k+2 , k+3 , \dots , 2k - 1 \}
\]
is a s.c.sf. set in $\mathbb{Z}_{3k+1}$ with $|S| = k$ (see
\cite{Haviv-Levy}).  In the case that 
$k = 3$, we have $n = 10$ and the Petersen graph is 
a 10-vertex 3-regular $K_3$-saturated graph.  
Finally, suppose 
that $n = 3k$.  For $k \geq 7$, we consider $\mathbb{Z}_{3k}$ and the 
subset 
\[
S = \{ k -1 , k+1 , 2k + 1 , 2k-1 \} \cup
\{ k+3 , k+4 ,k+5, \dots , 2k-3 \}.
\]
By Theorem 3.7 of \cite{Haviv-Levy} (with $n = 3k$, $s = k -1$, and $T = \{0,2 \}$ in the notation from that paper), $S$ is a 
s.c.sf. 
set with $|S| = k-1$.  Thus, $\textup{Cay}( \mathbb{Z}_{3k} , S )$ will be a $K_3$-saturated graph 
with $3k$ vertices that is $(k-1)$-regular.  
\end{proof}

\bigskip
Write $G(k,r)$ for the graph constructed in Theorem 1.1 in the 
case that $n = 3k+r$ where (i) $k \geq 1$ if $r =2$, (ii) $k \geq 4$ 
if $r = 1$, and (iii) $k \geq 7$ if $r =0$.  If $k = 3$ and $r = 1$, 
let $G(3,1)$ denote the Petersen graph.  In all cases, if $k  \geq 7$, then $G(k,r)$ is a 
\begin{center}
$(3k+r)$-vertex $(k+(r-1))$-regular $K_3$-saturated graph for $r \in \{0,1,2 \}$.
\end{center}

The graphs $G(k,r)$ provide us with a family of graphs for which we can apply 
Proposition \ref{useful tool}.  Our approach also needs a special $K_4$-saturated graph.  

\begin{proposition}\label{5k+4}
Let $k \geq 3$ be an integer and $n = 5k+4$.  If 
\[
S= \{1,2,n-1,n-2 \} \cup \left\{ x : 6 \leq x \leq 5(k-1)+3 
~ \mbox{and} ~ x \equiv 1 ~\mbox{or}~~3 ( \textup{mod}~5) \right\},
\]
then the graph $H(k) : = \textup{Cay}( \mathbb{Z}_n , S )$ is a $(2k+2)$-regular
$K_4$-saturated graph $n$ vertices.  
\end{proposition}

Before giving the proof of Proposition \ref{5k+4}, we give an example
for the readers' convenience.  
If $k = 6$, then $n = 34$ and 
\[
S = \{1,2,33,32 \} \cup \{ 6,8,11,13,16,18,21,23,26,28 \}.
\]

\smallskip

\begin{proof}[Proof of Proposition \ref{5k+4}]
Let $N_0$ be the neighborhood of $0$ in $H(k)$ so that $N_0 = S$.
Let $H_0 (k)$ be the subgraph of $H(k)$ induced by $S$.  
Since $H(k)$ is vertex transitive, to show that $H(k)$ is 
$K_4$-free it is enough to prove that $H_0(k)$ is $K_3$-free.
For the saturation property, it is enough to show that 
for any vertex $x$ not adjacent to 0, there is an edge $\{y,z \}$ in 
$H_0(k)$ where both $y$ and $z$ are adjacent to $x$ in $H(k)$.  

In $H_0(k)$, the vertex 1 is adjacent to only 2 and $n-1$.
Since $(n -1) -2  = n -3 \notin S$, 2 is not adjacent to $n-1$ so there is no triangle in 
$H_0(k)$ using vertex 1. By symmetry, there is 
no triangle in $H_0(k)$ containing $n - 1$.  No three vertices 
from the set 
\[
S':= \{ x : 6 \leq x \leq 5 (k-1) + 3
~ \mbox{and} ~ x \equiv 1 ~\mbox{or}~3 ( \textup{mod}~5) \}
\]
can form a triangle since given any three distinct vertices $x_1, x_2, x_3 \in S'$, 
at least one of the differences $x_1 - x_2$, $x_2 - x_3$, $x_3 - x_1$ will
be divisible by 5.
If say $x_1  - x_2$ is divisible by 5, then $x_1$ and $x_2$ are not adjacent since $S$ contains no element divisible by 5.  

Vertices 2 and $n-2$ are not adjacent since $4$ is not in $S$, so no 
triangle in $H_0 (k)$ can contain both 2 and $n-2$.  The only remaining case
to check is if 2 (or $n-2$) and two vertices $x_1,x_2 \in S'$ form a 
triangle.  Assuming that $x_1$ and $x_2$ are adjacent and both are in $S'$, then without loss of 
we can say that $x_1 = 5 \ell +1$.  Note that 
$x_1 - 2 = 5 ( \ell - 1) +4$ and $S$ does not contain any element 
congruent to 4 modulo 5.  Hence, 2 is not adjacent to $x_1$.  We 
conclude that $H_0 (k)$ is $K_3$-free and so $H(k)$ is $K_4$-free.

Now we show that given any vertex $x $ not adjacent to 0, 
$x$ is in a triangle $\{x,y,z \}$ where $y$ and $z$ are vertices in $H_0 (k)$.  
For each $x$, we give such a triangle. In some cases we compute differences to verify an adjacency.  
\begin{itemize}
\item $x = 3$: A triangle containing $3$ is $\{1,2,3 \}$.  
Note the differences $3 - 2 = 1$, $3 - 1 = 2$, $2 - 1=1$ are all in $S$ so that $\{1,2,3 \}$ is a triangle in $H(k)$.  
\item $x = 4$: A triangle containing 4 is $\{ 6 , n - 2 , 4 \}$.  
Note $6  - 4 = 2 \in S$, 
\begin{center}
$(n -2) - 4 = 2k - 2 = 5 ( k - 1) + 3 \in S$ ~~~ and ~~~
$(n -2) - 6 = 5(k-1)+1 \in S$.
\end{center} 
\item $x = 5$: A triangle containing 5 is $\{ 5,11,13 \}$.  Note that for all $k \geq 3$, $\{2,6,8 \} \subseteq S$ so that the differences $13 - 5, 11-5,13-11$ are in $S$.  
\item The case $x \in \{ n - 3,n -4 , n - 5 \}$ follows from symmetry and the above three cases.
\end{itemize}
This covers all vertices not adjacent to 0 with the exception of all 
$x$ for which 
\begin{center}
$6 \leq x \leq 5 (  k - 1) + 3$ ~~and~~ $x \equiv 0,2, ~\mbox{or}~4 ( \textup{mod}~5)$.
\end{center}
First suppose $x \equiv 2 ( \textup{mod}~5)$, say $x = 5 \ell + 2$ with 
$1 \leq \ell \leq  k - 2$.  Then 
\[
\{ x , y =  5 \ell + 1 , z = 5 \ell + 3 \}
\]
is a triangle containing $x$ with $y,z \in H_0(k)$.  Indeed, the differences 
$x - y = 1$, $z - x = 1$, and $z - y = 2$ are all in $S$.  

Now assume that $x = 5 \ell + 4$ with $1 \leq \ell \leq k - 2$.  
Consider 
\[
\{ x = 5 \ell + 4 , y = 5 ( \ell +1 ) + 1 , z = n - 2 \}.
\]
We have $y -x = 2 \in S$, 
\[
z - x = 5k+2 - ( 5 \ell + 4) = 5 ( k - ( \ell + 1) ) + 3 \in S
\]
and
\[
z - y = 5k+2 - (5 \ell + 6) = 5 ( k - ( \ell + 1) ) + 1 \in S.
\]
Therefore, $\{ x,y,z \}$ is a triangle with $y,z \in H_0 (k)$.

Lastly, assume $x = 5 \ell + 0$ with $2 \leq \ell \leq k - 1$.  Then $\{ x = 5 \ell , y = 2 , z = 5 \ell  - 2\}$ is a triangle with $y,z \in S$.  The computations 
\begin{center}
$x - y = 5 ( \ell - 1) + 3 \in S$, ~~ $x - z = 2 \in S$, ~~~
$z - y = 5 ( \ell - 1) + 1 \in S$
\end{center}
show that $\{x,y,z \}$ is a triangle in $H (k)$ and by definition, $y,z \in S$.    This completes
the proof of Proposition \ref{5k+4}.  
\end{proof} 

\bigskip

Using Proposition \ref{5k+4} and $K_3$-saturated graphs $G(k,r$),
we can prove the main results of this section.

\begin{theorem}\label{all K4}
For all $n \geq 36$, there is a regular $K_4$-saturated graph on $n$ vertices.
\end{theorem}
\begin{proof}
Let $n \geq 36$.  
Proposition \ref{useful tool} implies that if $G$ is an $N$-vertex 
$K_3$-saturated $d$-regular graph,
then $G + \overline{K_{N-d}}$ is a $K_4$-saturated graph with 
$2N - d$ vertices that is $N$-regular.  
If we apply this with $G = G(k,r)$, then we obtain a $(3k+r)$-regular
$K_4$-saturated graph with 
\[
2( 3k+r) - ( k + ( r- 1) ) = 5k + r + 1 
\]
vertices. When $r = 0$, $G(k,r)$ exists for 
all $k \geq 7$. When $r = 1$, $G(k,r)$ exists 
for all $k \geq 4$ and we may use the Petersen 
graph to deal with the case of $k = 3$.  Lastly,
when $r = 2$, $G(k,r)$ exists for all $k \geq 1$.  
This implies that for any $n > 36$  with $n = 5k+r+1$ for some $r \in \{0,1,2 \}$, 
there is a $K_4$-saturated graph with 
$5k+r+1$ vertices that is $(3k+r)$-regular.

If $n \equiv 4 ( \textup{mod}~5)$, then
by Proposition \ref{5k+4} there is a $K_4$-saturated graph with 
$5k+4$ vertices that is $(2k+2)$-regular provided $k \geq 3$.  

The last possibility is when $n$ is divisible by 5.  Using Mathematica
\cite{Mathematica}
a computer search was used to construct regular $K_4$-saturated graphs 
with $n$ vertices for the following values of $n$: 
\begin{center}
$9, 11, 13,14,15,16,17,19,20,23,25,29,31,49$
\end{center}
Note that with this list together with the constructions given so far, 
there is a regular $K_4$-saturated graph with $p$ vertices for any prime $p \geq 11$, 
and also on 9, 16, 25, and 49 vertices.   

First suppose that 25 divides $n$, say $n = 25 t$.  Let $H$ be a $K_4$-saturated
regular graph with 25 vertices.  Let 
$H'$ be obtained from $H$ by replacing each vertex in $H$ with an independent set of size $t$.  We call this the \emph{$t$-blow up of $H$}.  The graph $H'$ is a
$K_4$-saturated graph with $25t$ vertices.  The method of blowing up vertices into 
independent sets has been used before (see \cite{Gerbner-et-al}).  

Now suppose that 5 divides $n$, but 25 does not divide $n$.  Let $n = 5s$.  
If $s$ is divisible by some prime $p \geq 11$ or by 9, 16, or 49, 
then we can use a $t$-blow up to construct a regular $K_4$-saturated graph with $n$ vertices.  
Hence, we can assume that $n =  2^a 3^b 5^1 7^c$ where
$0 \leq a \leq 3$, $0 \leq b , c \leq 1$.  This tells that 
$n$ is a divisor of 840.  All divisors of $2^3 3^1 5^1 7^1 = 840$ that are greater than 
35 are divisible by either 14, 15, or 20.  This implies that we may again
complete the proof using $t$-blow ups of regular $K_4$-saturated graphs on 14, 15, or 20 vertices.  
\end{proof}

\smallskip

Now we prove our result for $K_5$-saturated graphs.  

\begin{theorem}\label{k5 sat}
For all $n  \geq 59$, there is a regular $K_5$-saturated graph with $n$ vertices.
\end{theorem}
\begin{proof}
Fix an integer $k \geq 7$ and integers 
$r_1 ,r_2 \in \{0,1,2 \}$.  Consider the 
join $G( k , r_1)  + G(k,r_2)$ of $G(k,r_1)$ and 
$G(k, r_2)$. This graph will be 
$K_5$-saturated, have $6k + r_1 + r_2$ vertices and since
\[
(3k + r_1) + (k + ( r_2 - 1) ) = (3k+r_2) + ( k + ( r_1 - 1) ),
\] 
the graph $G(k,r_1) + G(k,r_2)$ will be $d : =  4k + r_1 + r_2 - 1$-regular.  
This shows that for any integer $k \geq 7$ and 
$r_1,r_2 \in \{0,1,2 \}$, there is a regular $K_5$-saturated  
graph with $6k + r_1 + r_2$ vertices.  This covers all integers
$n \geq 42$ with the exception of those $n$ for which $n \equiv 5 ( \textup{mod}~6)$.
For this case, we use the construction of a s.c.sf. set 
given in \cite{Haviv-Levy}.  For $k \geq  9$, in $\mathbb{Z}_{3k+2}$ let 
\[
S' = \{ k - 1 , k +2 , k +3 \} 
\cup \{ k+5,k+6,k+7,\dots , 2k-3 \} 
\cup \{ 2k-1 , 2k+1 , 2k+3 \}.
\]
If $G'(k) = \textup{Cay}( \mathbb{Z}_{3k+2} , S')$, then $G'(k)$ is 
a $K_3$-saturated graph with $3k+2$ vertices that is $(k-1)$-regular.
Thus, for $k \geq 9$, 
\[
G'(k) + G(k+1,0)
\]
is a $K_5$-saturated graph with $3k+2 + 3(k+1)= 6k+5$ vertices.
Since
\[
(3k+2) + k  = 3(k+1) + k - 1,
\]
this graph is regular of degree $(k-1) + 3(k+1) = 4k+2$.  Hence, for 
any $n \geq 59$ with $n \equiv 5 ( \textup{mod}~6)$, there 
is a regular $K_5$-saturated graph.  
\end{proof}

\smallskip

We conclude this section by using Proposition \ref{useful tool} and 
our constructions given so far to prove Theorem \ref{all large values}.

\begin{proof}[Proof of Theorem \ref{all large values}]
We will only prove part (i) in  because the proof of (ii) is similar.
We use induction on $\delta$ with the base case of $\delta = 1$
given by Theorem \ref{all K3}.    
For $r' \in \{0,1,2 \}$ and large enough $n$, we let $G_1$ be 
a $K_3$-saturated graph with $3k+r'$ vertices that is regular of 
degree $k + (r'-1)$.  For $r'' \in \{0,1,2 , \dots , 2 ( \delta -1 ) \}$
and large enough $n$, let $G_2$ be a $K_{ 2 \delta - 1}$-saturated 
graph with $3 ( \delta - 1) k + r''$ vertices that 
is $[ (3 ( \delta - 1) - 2) k + (r'' - 1) ]$-regular. Such a graph 
exists by the induction hypothesis.   
Applying Proposition \ref{useful tool} to $G_1$ and $G_2$ gives a 
$K_{2 \delta +1}$-saturated graph with $3 \delta k + r' +r''$ vertices 
that is $[ ( 3 \delta - 2 )k + (r' + r'')-1]$-regular.  The sum 
$r' + r''$ ranges over all values in $\{0,1,2, \dots , 2 \delta \}$ as $r'$ and $r''$ range
over their possible values.  

\smallskip
For proving part (ii), note that our results on $K_4$-saturation prove 
the base case $\delta = 1$ (see the proof of Theorem \ref{all K4}).  
The remainder of the proof is almost 
 identical to the proof of part (i).    
 \end{proof}


\section{Concluding Remarks}\label{conclusion}


The methods used in this paper relied on first finding sets $S \subseteq \mathbb{Z}_n$ with 
special properties.  All upper bounds obtained in this paper are quadratic in the number of vertices, 
with the exception of $\textup{rsat}(n , K_3) = \Theta (n^{3/2})$ which follows from 
\cite{Gerbner-et-al} and \cite{Haviv-Levy}.  Finding subquadratic upper bounds, 
assuming such bounds exist, would hopefully 
lead to other interesting constructions in groups.    

A close look at the proof of Theorem \ref{super sum set} shows $(2 \alpha + 2)S = \mathbb{Z}_n \backslash S$ and
that 
\[
( (2 \alpha + 2) S) \backslash \{ 0 \} = \mathcal{R}_{ 2 \alpha + 2} (S).  
\]
This motivates the following problem already mentioned in the Introduction: determine for which $n$ and $k$ is there a symmetric set $S \subseteq \mathbb{Z}_n$ such that 
\[
k S = \mathbb{Z}_n \backslash S.
\]
In the literature, such sets are called \emph{complete}.  In this paper, we gave a family of examples in Theorem 
\ref{super sum set}.  
  Using 
the computer program Mathematica \cite{Mathematica}, we 
found more sets satisfying the assumptions of Proposition \ref{construction} 
for $k =4$ and these give rise to regular $C_5$-saturated graphs.  All of these examples are 
in $\mathbb{Z}_n$.  

\begin{center}
\begin{tabular}{c | c || c | c  }
$n$ & $S$ & $n$ & $S$ \\ \hline
17 & $\{ 1,3,14,16 \}$ &   37 & $\{1,3,5,7,30,32,34,36 \}$ \\
21 & $\{ 1, 6, 8, 13,15,20 \}$  & 39 & $\{1,3,14,25,36,38 \}$  \\
23 & $\{ 1,5,18,22 \}$ & 41 & $\{1,5,11,30,36,40 \}$ \\
25 & $\{ 1,7,18,24 \}$ & 43 & $\{1,6,8,35,37,42 \}$ \\
27 & $\{1,3,5,22,24,26 \}$ & 45 & $\{1,6,8,37,39,44 \}$ \\
29 & $\{1 , 12 ,17,28 \}$ & 47 & $\{1,3,13,34,44,46 \}$   \\
33 & $\{ 1,3,7,26,30,32 \}$ & 49 & $\{1,3,19,30,46,48 \}$ \\
37 & $\{1,3,5,7,30,32,34,36 \}$   & 51 & $\{1,12,23,28,39,50 \}$  \\
39 & $\{1,3,14,25,36,38 \}$ &~ & ~ \\
\end{tabular}
\bigskip

Sets $S \subset \mathbb{Z}_n$ the generate an $|S|$-regular $C_5$-saturated $n$-vertex graph  
\end{center} 

Since a complete bipartite graph with part sizes at least 3 is $C_5$-saturated,
$\textup{rsat}(n, C_5)$ exists for all even $n$.  The table above 
shows that $\textup{rsat}(n , C_5)$ exists for all 
odd $n$ from 17 to 51 with the exception of 19, 31, and 35 where our 
computer search failed to find such sets.  While it is arguable how 
compelling this evidence is, we are very tempted to conjecture 
that one can find a set $S \subset \mathbb{Z}_n$ satisfying 
the requirements of Proposition \ref{construction} for $k=4$ and for all sufficiently large 
odd $n$.  Of course if true, this would 
imply $\textup{rsat}(n , C_5)$ exists for all large enough $n$.
Sporadic examples
satisfying Proposition \ref{construction} for $k= 5$ were also found, but we were not able 
to find a pattern.  This problem certainly seems interesting on its own. 

\begin{question}\label{question}
 For which even integers $k \geq 4$ and $n$ is 
there a symmetric 
subset $S \subset \mathbb{Z}_n$ such that $0 \notin (k+1)S$ and 
\[
\mathcal{R}_k(S)  = \mathbb{Z}_n \backslash (S \cup \{0 \} ).
\]
\end{question}

For large enough $n$, $\mathbb{Z}_n$ contains 
a complete sum-free set $S \subseteq \mathbb{Z}_n$ (see \cite{Haviv-Levy}).
Such a set
will satisfy $0 \notin 3S$ and 
$\mathcal{R}_2 ( S) = \mathbb{Z}_n \backslash ( S \cup \{ 0 \} )$
so when $k = 2$, the answer to Question \ref{question} is 
``all sufficiently large $n$".



\end{document}